\documentclass{nyj}
\usepackage{url}
\title{SLE Coordinate Changes}
\author{Oded Schramm}
\author{David B. Wilson}
\begin{document}
\begin{abstract}
  The purpose of this note is to describe a framework which unifies
  radial, chordal and dipolar SLE.  When the definition of
  SLE$(\kappa;\rho)$ is extended to the setting where the force points
  can be in the interior of the domain, radial SLE$(\kappa)$ becomes
  chordal SLE$(\kappa;\rho)$, with $\rho=\kappa-6$, and vice versa.
  We also write down the martingales describing the Radon--Nykodim
  derivative of SLE$(\kappa;\rho_1,\dots,\rho_n)$ with respect to
  SLE$(\kappa)$.
\end{abstract}
\maketitle
\tableofcontents

\newtheorem{theorem}{Theorem}
\newtheorem{conj}{Conjecture}
\newtheorem{lemma}[theorem]{Lemma}
\newtheorem{proposition}[theorem]{Proposition}
\newtheorem{corollary}[theorem]{Corollary}
\theoremstyle{definition}\newtheorem{remark}[theorem]{Remark}
\theoremstyle{definition}\newtheorem{definition}[theorem]{Definition}
\newcommand{\Z}{\mathbb{Z}}
\newcommand{\R}{\mathbb{R}}
\newcommand{\C}{\mathbb{C}}
\renewcommand{\H}{\mathbb{H}}
\newcommand{\N}{\mathbb{N}}
\newcommand{\E}{\mathbb{E}}
\newcommand{\D}{\mathbb{D}}
\newcommand{\eps}{\varepsilon}
\newcommand{\Var}{\operatorname{Var}}
\newcommand{\Cov}{\operatorname{Cov}}
\newcommand{\SLE}{\operatorname{SLE}}
\newcommand{\Rad}{\operatorname{Rad}}
\renewcommand{\Cap}{\operatorname{Cap}}
\newcommand{\rad}{\operatorname{rad}}
\newcommand{\lcap}{\operatorname{cap}}
\renewcommand{\mod}{\operatorname{mod}}
\newcommand{\wind}{\operatorname{wind}}
\newcommand{\twist}{\operatorname{twist}}
\newcommand{\note}[1]{\begin{center}\fbox{\parbox{5in}{\textbf{Note: #1}}}\end{center}}
\newcommand{\old}[1]{}
\newcommand{\aref}[1]{Appendix~\ref{sec:#1}}
\newcommand{\sref}[1]{\S~\ref{sec:#1}}
\newcommand{\tref}[1]{Theorem~\ref{thm:#1}}
\newcommand{\lref}[1]{Lemma~\ref{lem:#1}}
\newcommand{\pref}[1]{Proposition~\ref{pro:#1}}
\newcommand{\cref}[1]{Corollary~\ref{cor:#1}}
\newcommand{\fref}[1]{Figure~\ref{fig:#1}}
\newcommand{\eref}[1]{Equation~\eqref{eqn:#1}}
\newcommand{\tbref}[1]{Table~\ref{tbl:#1}}
\newcommand{\W}{\mathring{W}}
\newcommand{\z}{\mathring{z}}
\newcommand{\g}{\mathring{g}}
\newcommand{\mob}{M\"obius}
\newcommand{\ito}{It\^o}
\newcommand{\SLErad}{\operatorname{SLE}^{\operatorname{rad}}}
\newcommand{\SLEk}{SLE{$_\kappa$}}
\newcommand{\p}{\partial}
\newcommand{\closure}[1]{\overline{#1}}
\newcommand{\dist}{\operatorname{dist}}
\newcommand{\KK}{\mathcal{G}}
\newcommand{\QED}{\qed\medskip}
\newcommand{\normalizedmap}{\Phi}
\renewcommand{\P}{{\bf P}}
\newcommand{\md}{\mid}
\newcommand{\Bb}[2]{{\renewcommand{\md}{\bigm| }#1\bigl[#2\bigr]}}
\newcommand{\BB}[2]{{\renewcommand{\md}{\Bigm| }#1\Bigl[#2\Bigr]}}
\newcommand{\Bs}[2]{{\renewcommand{\md}{\mid}#1[#2]}}
\newcommand{\Pb}{\Bb\P}
\newcommand{\Eb}{\Bb\E}
\newcommand{\PB}{\BB\P}
\newcommand{\EB}{\BB\E}
\newcommand{\Ps}{\Bs\P}
\newcommand{\Es}{\Bs\E}

\newcommand{\RSslecont}{math.PR/0106036}
\newcommand{\LSWrest}{MR1992830}
\newcommand{\LSWi}{MR2002m:60159a}
\newcommand{\LSWii}{MR2002m:60159b}
\newcommand{\LSWiii}{MR1899232}
\newcommand{\LSWlesl}{MR2044671}
\newcommand{\RevuzYor}{MR2000h:60050}
\newcommand{\DubedatMartingalesDuality}{MR2118865}
\newcommand{\WernerTransformations}{MR2060031}
\newcommand{\LWloops}{MR2045953} 
\newcommand{\diPolarSLE}{MR2140124}
\newcommand{\zigzagSLE}{math-ph/0401019}
\newcommand{\DubedatCommutation}{math.PR/0411299}
\newcommand{\Langetal}{MR1230963}
\newcommand{\SmirnovPerc}{MR1851632}
\newcommand{\WernerSurvey}{MR1905353}
\newcommand{\CardySurvey}{MR2037564}
\newcommand{\CardySurveySLE}{MR2148644}
\newcommand{\SWpercexpo}{MR1879816}
\newcommand{\SchSLE}{MR1776084}
\newcommand{\CardyFormula}{MR92m:82048}
\newcommand{\KenyonLERW}{MR1819995}
\newcommand{\WernerStFlour}{MR2079672}
\newcommand{\LawlerSLEintro}{LawlerSLEintro}
\newcommand{\LawlerSLE}{LawlerSLE}
\newcommand{\KagerNienhuis}{MR2065722}
\newcommand{\SLEsurveys}{\WernerStFlour,\LawlerSLE,\KagerNienhuis,\CardySurveySLE}

\section{Introduction}
SLE, or Stochastic Loewner Evolution~\cite{\SchSLE}, describes random paths
in the plane by specifying a differential equation satisfied by the
conformal maps into the complement.
The main interest in these paths stems from their relationship with
scaling limits of critical models from statistical physics and two
dimensional Brownian motion.
For motivation and background on SLE, the reader is advised to consult
the surveys~\cite{\SLEsurveys} and the references cited there.

SLE comes in several different flavors, and the principal goal of this
short paper is to observe a unification of these different variations.
Basically, SLE is Loewner evolution with Brownian motion as its driving
paramenter.
More explicitly, chordal SLE is defined as follows:
fix $\kappa\ge 0$.
Let $B_t$ be one-dimensional standard Brownian motion, started as
$B_0=0$.  Set $W_t=\sqrt\kappa\, B_t$.
If we fix any $z\in\H$, we may consider the solution of
the ODE
\begin{equation}
\label{chordal}
\p_t g_t(z) = \frac 2{g_t(z)-W_t},\qquad g_0(z)=z.
\end{equation}
Let $\tau=\tau_z$ be the supremum of the
set of $t$ such that $g_t(z)$ is well-defined.
Then either $\tau_z=\infty$, or $\lim_{t\nearrow\tau} g_t(z)-W_t=0$.
Set $H_t:=\{z\in\H:\tau_z>t\}$.
It is easy to check that for every $t>0$ the map $g_t:H_t \to \H$ is conformal.
It has been shown~\cite{\RSslecont,\LSWlesl}
 that the limit $\gamma(t)=\lim_{z\to W_t}g_t^{-1}(z)$ exists and
is continuous in $t$.  This is the path defined by SLE.
It is easy to check that for every $t\ge0$, $H_t$ is the unbounded connected
component of $\H\setminus \gamma[0,t]$.


Radial SLE is defined in a similar manner, except that the upper half
plane $\H$ is replaced by the unit disk $\D$, and the differential
equation~\eqref{chordal} is replaced by
\begin{equation}
\label{radial}
\p_t g_t(z) = -g_t(z)\,\frac {g_t(z)+W_t}{g_t(z)-W_t},\qquad g_0(z)=z,
\end{equation}
where $W_t$ is Brownian motion on the unit circle $\p \D$, with time
scaled by a factor of $\kappa$.
The definition is presented in detail and generalized in Section~\ref{s.unified}.
Yet another version of SLE, dipolar SLE, was introduced in~\cite{\zigzagSLE}
and further studied in~\cite{\diPolarSLE}.

In \cite{\LSWii} it was shown that radial SLE and chordal SLE have equivalent
(i.e., mutually absolutely continuous) laws, up to a time change, when stopped before
the disconnection time (the precise meaning of this should become clear soon).
Here, we carry this calculation a step further, and describe what a radial SLE looks
like when transformed to the chordal coordinate system and vice versa.
We extend the
definition of the SLE$(\kappa;\rho)$ processes from~\cite{\LSWrest,\DubedatMartingalesDuality}
to the setting where the force points are permitted to
be in the interior of the domain.
With the extended definitions, the SLE$(\kappa;\rho)$
radial processes are transformed into appropriate SLE$(\kappa;\rho)$
chordal processes, and vice versa.

\medskip

We start with the definition of chordal SLE$(\kappa;\rho)$ with a force point in
the interior.

\begin{definition}\label{d.chordali}
Let $z_0\in\H$, $w_0\in\R$, $\kappa\ge 0$, $\rho\in\R$.
Let $B_t$ be standard one-dimensional Brownian motion.  Define $(W_t,V_t)$ to be the
solution of the system of SDEs
\begin{align*}
d W_t &= \sqrt\kappa \,dB_t+\rho\,\Re \frac 1{W_t-V_t}\,dt,\\
d V_t &= \frac {2\,dt}{V_t-W_t},
\end{align*}
starting at $(W_0,V_0)=(w_0,z_0)$, up to the first time $\tau>0$
such that \[\inf\{|W_t-V_t|:t\in[0,\tau)\}=0.\]
(Here, $\Re z$ denotes the real part of $z$.)
Then the solution of Loewner's chordal equation~\eqref{chordal}
with $W_t$ as the driving term will be called chordal SLE$(\kappa;\rho)$
starting at $(w_0,z_0)$.
\end{definition}

One motivation for this definition comes from the fact that, as we will later
see, ordinary radial SLE$(\kappa)$ starting at $W_0=w$ is transformed by
a \mob\ transformation $\psi:\D\to\H$ to a time changed chordal SLE$(\kappa;\kappa-6)$
starting at $(\psi(w),\psi(0))$ (both up to a corresponding positive stopping
time).

In~\cite{\WernerTransformations} it was shown that a chordal SLE$(\kappa;\rho)$ process
can be viewed as an ordinary SLE$(\kappa)$ process weighted by a martingale.  In
Section~\ref{s.mart}
we extend this to SLEs with force points in the interior of the domain.

\section{A unified framework for SLE$(\kappa;\rho)$}
\label{s.unified}

In order to produce a clean and general formulation of the change of coordinate results, we
need to introduce a bit of notation which will enable us to deal with the
chordal and radial versions simultaneously.
Let $X\in\{\D,\H\}$.  We let $\Psi_X(w,z)$ denote the corresponding Loewner vector field;
that is,
$$
\Psi_{\D}(w,z)= -z\,\frac {z+w}{z-w},\qquad
\Psi_{\H}(w,z)= \frac {2}{z-w}.
$$
Thus, Loewner's radial [respectively chordal] equation may be written as
$$
\p_t g_t(z)=\Psi_X\bigl(W_t,g_t(z)\bigr)
$$
if $X=\D$ [respectively $X=\H$].
Let $I_\D$ denote the inversion in the unit circle $\p\D$, and let
$I_\H$ denote the inversion in the real line $\p\H$, i.e.,
complex conjugation.
Set
$$
\tilde\Psi_X(z,w):=\frac{\Psi_X(z,w)+\Psi_X(I_X(z),w)}2.
$$
In radial ordinary SLE, $W_t=W_0\exp(i\sqrt\kappa B_t)$, where $B_t$ is a Brownian motion.
By \ito's formula, $dW_t=-(\kappa/2)\, W_t\,dt+i\,\sqrt\kappa\,W_t\,dB_t$.
On the other hand, in ordinary chordal SLE, one simply has $dW_t=\sqrt\kappa\,dB_t$.
We thus set
$$
\KK_\D(W_t,dB_t,dt):=
-(\kappa/2)\, W_t\,dt+i\,\sqrt\kappa\,W_t\,dB_t,\qquad
\KK_\H(W_t,dB_t,dt):=\sqrt\kappa\,dB_t.
$$

\begin{definition}\label{d.genslerho}
Let $\kappa\ge 0$, $m\in\N$, $\rho_1,\rho_2,\dots,\rho_m\in\R$.
Let $X\in\{\H,\D\}$ and let $V^1,V^2,\dots,V^m\in\overline {X}$.
(When $X=\H$, we include $\infty$ in $\overline{X}$.)
Let $w_0\in \partial X\setminus\{\infty,V^1,\dots,V^m\}$
and let $B_t$ be standard one-dimensional Brownian motion.
Consider the solution of the SDE system
\begin{align}
dW_t&=\KK_X(W_t,dB_t,dt)+\sum_{j=1}^m\frac{\rho_j}2\,\tilde\Psi_X(V^j_t,W_t)\,dt,
\label{e.dW}
\\
dV^j_t&=\Psi_X(W_t,V^j_t)\,dt,\qquad j=1,2,\dots,m
\end{align}
starting at $W_0=w_0$ and $V^j_0=V^j$, $j=1,\dots,m$,
up to the first time $\tau$ such that for some $j$
$\inf\{|W_t-V^j_t|:t<\tau\}=0$.
Let $g_t(z)$ be
the solution of the ODE
$$
\p_t g_t(z)=\Psi_X(W_t,g_t(z))
$$
starting at $g_0(z)=z$.  Let $K_t=\{z\in\closure X:\tau_z\le t\}$,
$t<\tau$, be the corresponding hull,
defined in the same way as for ordinary SLE.  Then the evolution
$t\mapsto K_t$ will be called $X$-SLE$(\kappa;\rho_1,\dots,\rho_m)$
starting at $(w_0,V^1,\dots,V^m)$.  When $X=\D$, we refer to $\D$-SLE
as radial SLE, while $\H$-SLE is chordal.
The points $V^j_t$ will be called {\bf force points}.
\end{definition}

Note that $V_t^j=g_t(V^j)$, for all $t\in[0,\tau)$.

In some situations, it is possible and worthwhile to extend the definition
of SLE$(\kappa;\rho_1,\dots,\rho_m)$ beyond the time $\tau$
(see, e.g.,~\cite{\LSWrest}), but we do not deal with this here.

\section{Some \mob\ coordinate changes}

Set $o_\H:=\infty$ and $o_\D:=0$ and as before $X\in\{\H,\D\}$.
Observe that if $V^m=o_X$,
 then the value of $\rho_m$ has
no effect on the $X$-{\rm SLE}$(\kappa;\rho_1,\dots,\rho_m)$,
and, in fact, the {\rm SLE}$(\kappa;\rho_1,\dots,\rho_m)$
reduces to SLE$(\kappa;\rho_1,\dots,\rho_{m-1})$.
Consequently, by increasing $m$ and appending $o_X$ to the force points, if necessary,
there is no loss of generality in assuming that
$\sum_{j=1}^m\rho_j=\kappa-6$.  This assumption simplifies somewhat the
statement of the following theorem:

\begin{theorem}\label{t.coord}
Let $X,Y\in\{\D,\H\}$, and let $\psi:X\to Y$ be a conformal homeomorphism;
that is, a \mob\ transformation satisfying $\psi(X)=Y$.
Let $w_0\in\partial{X}\setminus\{\infty\}$, $V^1,V^2,\dots,V^m\in
\overline {X}\setminus\{w_0\}$.
Suppose that $\rho_1,\rho_2,\dots,\rho_m\in\R$ satisfy $\sum_{j=1}^m\rho_j=\kappa-6$.
Then the image under $\psi$ of the $X$-SLE$(\kappa;\rho_1,\dots,\rho_m)$
starting from $(w_0,V^1,\dots,V^m)$ and stopped at some a.s.\ positive stopping
time has the same law
as a time change of the $Y$-{\rm SLE}$(\kappa;\rho_1,\dots,\rho_m)$ starting from
$(\psi(w_0),\psi(V^1),\dots,\psi(V^m))$
stopped at an a.s.\ positive stopping time.
\end{theorem}

The stopping time for the $X$-SLE may be taken as the minimum
of $\tau$ and the first time $t$ such that $\psi^{-1}(o_Y)\in K_t$.

The case when $X=Y=\H$ and the force points are on the boundary appears in \cite{\DubedatCommutation}.

\proof
Suppose first that $X=\H$ and $Y=\D$.
Let $g_t$ be the family of maps $g_t:\H\setminus K_t\to\H$.
Let $\phi_t:\H\to\D$ be the \mob\ transformation
such that $F_t:= \phi_t\circ g_t\circ \psi^{-1}$
satisfies the radial normalization $F_t(0)=0$ and
$F'_t(0)\in(0,\infty)$.
Let $z_t=x_t+i\,y_t:=g_t\circ\psi^{-1}(0)$, with $x_t,y_t\in\R$.
There is some (unique) $\lambda_t\in\p\D$ such that
$$\phi_t(z)=\lambda_t \frac{z-z_t}{z-\overline z_t}.$$
Let $s(t):=\log F'_t(0)$ and let $t(s)$ denote the inverse
of the map $t\mapsto s(t)$.
By the chain rule,
\begin{align*}
s'(t) = \p_t \log F'_t(0)
&=\p_t\log\biggl(\frac{\lambda_t}{2 i y_t}\,g_t'(\psi^{-1}(0)) \, (\psi^{-1})'(0)\biggr)\\
&=\Re \bigl(\p_t \log g_t'(\psi^{-1}(0))\bigr)-\p_t \log y_t.
\end{align*}
We have by~\eqref{chordal}
$$
\p_t y_t = \frac{-2 \,y_t}{|z_t-W_t|^2},
$$
 and
$$
\p_t g_t'(z)=\p_z \p_t g_t(z)= \p_z \frac{2}{g_t(z)-W_t}
= \frac {-2\, g_t'(z)}{(g_t(z)-W_t)^2}.
$$
Using this with $z=o$ in our previous expression for $s'(t)$ gives
$$
s'(t)=
\Re \frac{-2}{(z_t-W_t)^2} + \frac 2{|z_t-W_t|^2}
=
\frac  {4\,y_t^2} {|z_t-W_t|^4}.
$$

Set $\hat g_s(z):=F_{t(s)}(z)$,
$\hat W_s:= \phi_{t(s)}(W_{t(s)})$
and $\hat K_s:=\psi(K_{t(s)})$.
We need to verify that Loewner's equation holds:
\begin{equation}\label{e.hatL}
\p_s \hat g_s (z) =\Psi_\D\bigl(\hat W_s,\hat g_s(z)\bigr),\qquad z\in \D\setminus \hat K_s.
\end{equation}
This may be verified by brute-force calculation
or may be deduced from Loewner's theorem, as follows:
first, $\hat g_s$ is appropriately normalized: $\hat g_s(0)=0$ and $\hat g'_s(0)=e^s$.
Second, the domain of definition of $\hat g_s$ is clearly $\D\setminus\hat K_s$.
Lastly, if $\eps>0$ is small, then $g_{t+\eps}\circ g_t^{-1}$ is defined
on the complement in $\H$ of a small neighborhood of $W_t$.
Consequently, $\hat g_{s+\eps}\circ\hat g_s^{-1}$ is defined
on the complement in $\D$ of a small neighborhood of $\hat W_s$ when $\eps>0$
is small.  Therefore, Loewner's theorem gives~\eqref{e.hatL}.

It remains to calculate $d\hat W_s$.
Since
$
\hat W_s=\lambda_t\,\frac{W_t-z_t}{W_t-\overline z_t },
$
the \ito\ differential $d\log\hat W_s$ satisfies
\begin{equation}\label{e.hWs}
d\log\hat W_s=d\log\lambda_t+d\log({W_t-z_t})-d\log({W_t-\overline z_t }).
\end{equation}
Since $\hat g_s(\psi(\infty))=\phi_{t(s)}(\infty)=\lambda_{t(s)}$,
Equation~\eqref{e.hatL} with $z=\psi(\infty)$ gives
$
\p_s\lambda_{t(s)} = \Psi_\D(\hat W_s,\lambda_{t(s)})
$,
so
$$
A_0:=
d\log\lambda_t =
\lambda_t^{-1}\,\Psi_\D(\hat W_s,\lambda_{t})\,ds
=
-\hat W_s^{-1}\,\Psi_\D(\lambda_t,\hat W_s)\,ds.
$$
\ito's formula gives
\begin{equation}\label{e.f1}
\begin{aligned}
d\log(W_t-z_t)
&
=-\frac{\kappa\,dt}{2(W_t-z_t)^2}+\frac{dW_t-dz_t}{W_t-z_t}
\\&=
-\frac{\kappa\,dt}{2(W_t-z_t)^2}+\frac{dW_t}{W_t-z_t}+\frac{2\,dt}{(z_t-W_t)^2}
=
\frac{(4-\kappa)\,dt}{2(W_t-z_t)^2}+\frac{dW_t}{W_t-z_t}.
\end{aligned}
\end{equation}
Likewise,
\begin{equation}\label{e.f2}
-d\log(W_t-\overline z_t )
=
-\frac{(4-\kappa)\,dt}{2(W_t-\overline z_t )^2}-\frac{dW_t}{W_t-\overline z_t }.
\end{equation}
We now handle the sum of the $dt$
terms in~\eqref{e.f1} and~\eqref{e.f2}, using our
above expression for $s'(t)$.
\begin{equation*}
\begin{aligned}
A_1:=
\frac{(4-\kappa)\,dt}{2(W_t-z_t)^2}-
\frac{(4-\kappa)\,dt}{2(W_t-\overline z_t)^2}
&
= \frac{2\,i\,(4-\kappa)\,y_t\,(W_t-x_t)\,s'(t)^{-1}\,ds}
{|W_t-z_t|^4}
\\&
=
\frac{i\,(4-\kappa)\,(W_t-x_t)\,ds}{2\,y_t}
\\&
=
-\frac {4-\kappa}2\,\hat W_s^{-1}\,\Psi_\D(\lambda_t,\hat W_s)\,ds.
\end{aligned}
\end{equation*}
Next, we look at the $dW_t$ terms.
First, write~\eqref{e.dW} in slightly different form
$$
dW_t=\sqrt\kappa\,dB_t+\sum_{j=1}^{2m}\frac{\rho_j}4\,\Psi_\H(v_t^j,W_t)\,dt,
$$
where $v_t^j:=V_t^j$, $v_t^{j+m}:=I_\H( {v_t^{j}})$
and $\rho_{j+m}=\rho_j$ for $j=1,2,\dots,m$.
The sum of the $dW_t$ terms in~\eqref{e.f1} and~\eqref{e.f2} is
$$
\begin{aligned}
A_2
&
:=
\frac{dW_t}{W_t-z_t}- \frac{dW_t}{W_t-\overline z_t}
=
\frac{2\,i\,y_t}{|W_t-z_t|^2}\Bigl(\sqrt\kappa\,dB_t+
\sum_{j=1}^{2m}\frac{\rho_j}4\,\Psi_\H(v_t^j,W_t)\,dt\Bigr)
\\&\qquad\qquad
=
i\,\sqrt{\kappa\,s'(t)}\,dB_t
+i\,{s'(t)^{-1/2}}
\sum_{j=1}^{2m}\frac{\rho_j}4\,\Psi_\H(v_t^j,W_t)\,ds.
\end{aligned}
$$
Now define
$$
\hat B_s:=\int_0^{t(s)}\sqrt{s'(t)}\,dB_t.
$$
Then $\langle\hat B\rangle_{s_0}=\int_0^{t(s_0)}s'(t)\,dt=s_0$, and therefore
$s\mapsto\hat B_s$ is standard Brownian motion.
Also set $\hat v_t^j:=\phi_t(v_t^j)$, and note that $\phi_t(I_\H(z))=I_\D(\phi_t(z))$,
since $\phi_t:\H\to\D$ is a \mob\ transformation.
To further translate our expression for $A_2$ to the $\D$ coordinate system, we calculate
$$
i\,s'(t)^{-1/2}\Psi_\H(v_t^j,W_t)=
\frac{i\,|z_t-W_t|^2}{y_t(W_t-v_t^j)}
=
\hat W_s^{-1}\,\Psi_\D(\hat v_t^j,\hat W_s)
-\hat W_s^{-1}\, \Psi_\D(\lambda_t,\hat W_s)
.
$$
Consequently, since $\sum_{j=1}^m \rho_j=\kappa-6$,
$$
A_2=
i\,\sqrt{\kappa}\,d\hat B_s -\frac{\kappa-6}2\,\hat W_s^{-1}\,\Psi_\D(\lambda_t,\hat W_s)\,ds
+
\hat W_s^{-1}\sum_{j=1}^m \frac{\rho_j}2 \tilde\Psi_\D(\hat v_t^j,\hat W_s)\,ds.
$$
Settting $\hat w_s:=\log \hat W_s$,  we get
$$
d\hat w_s=A_0+A_1+A_2=
i\,\sqrt{\kappa}\,d\hat B_s
+
\hat W_s^{-1}\sum_{j=1}^m \frac{\rho_j}2 \tilde\Psi_\D(\hat v_t^j,\hat W_s)\,ds.
$$
Now, \ito's formula gives
\begin{multline*}
d\hat W_s=d\exp(\hat w_s)= \exp(\hat w_s)\,d\hat w_s+\frac 12\, \exp(\hat w_s)\,d\langle\hat w_s\rangle
\\
=i\,\sqrt{\kappa}\,\hat W_s\,d\hat B_s
+ \sum_{j=1}^m \frac{\rho_j}2 \tilde\Psi_\D(\hat v_t^j,\hat W_s)\,ds
-
\frac {\kappa\,\hat W_s}2 \,ds.
\end{multline*}
This completes the proof in the case $X=\H$ and $Y=\D$.

If $X=\D$ and $Y=\H$, we may just reverse the above equivalence.
To handle the case $X=Y=\H$, we may just write the
\mob\ trasformation $\psi:X\to Y$ as a composition
$\psi=\psi_2\circ\psi_1$, where $\psi_1:\H\to\D$ and $\psi_2:\D\to \H$,
and appeal to the above situtations.  The case $X=Y=\D$ is similar.
\QED

\begin{remark}\label{r.conceptual}
It is possible to come up with a somewhat more conceptual version of
some parts of this proof.
First, the time change $s'(t)$ is the rate at which the radial capacity
changes with respect to the chordal capacity.  This is known to be
$\phi_t'(W_t)^2$.  Second, we have $\hat W_s=\phi_t(W_t)$.
\ito's  formula gives
$$
d\hat W_s=(\p_t\phi_t)(W_t)\,dt+\phi_t'(W_t)\,dW_t+(\kappa/2)\,\phi_t''(W_t)\,dt.
$$
The terms $\phi_t'(W_t)\,\Psi_\H(V^j_t,W_t)$ that arise by expanding
$dW_t$ should be thought of as the image under $\phi_t$ of the vector
field $z\mapsto \Psi_\H(V^j_t,z)$, evaluated at $\hat W_s$.
The vector field $z\mapsto\Psi_\H(V^j_t,z)$ is the Loewner vector
field in $\H$, and it should be mapped under $\phi_t$ to a multiple of
the Loewner vector field in $\D$ plus some vector field which preserves
$\D$.
\end{remark}

\section{Strip SLE}

For comparison and illustration, we now mention another type of SLE.
It generalizes dipolar SLE as introduced in~\cite{\zigzagSLE}
and studied in~\cite{\diPolarSLE} as well as
a version of dipolar SLE with force points,
which was used in~\cite[\S3]{\LSWi}

\begin{definition}[Dipolar SLE]
Let $S=\{x+iy:x\in\R,\,y\in(0,\pi/2)\}$.
Set
\begin{align*}
\begin{matrix}
\Psi_S(w,z)=2\,\coth(z-w),\quad& I_S=I_\H, \\
\tilde\Psi_S(z,w)=\frac{\Psi_S(z,w)+\Psi_S(I_S(z),w)}2,\quad& \KK_S=\KK_\H.
\end{matrix}
\end{align*}
Then strip-SLE$(\kappa;\rho_1,\rho_2,\dots,\rho_m)$ is
defined using Definition~\ref{d.genslerho} with $X=S$
and $w_0\in\R$.
In the case where $m=0$, this coincides with dipolar SLE,
as defined in~\cite{\diPolarSLE}.
\end{definition}

Note that in strip-SLE$(\kappa;\rho_1,\dots,\rho_m)$,
possible force points at $+\infty$
and at $-\infty$ do exert an effect on the motion of $W_t$,
since $\Psi_S(\pm\infty,w)=\mp2$.
However, if $V^{m-1}=+\infty$ and $V^m=-\infty$, then
adding a constant to both $\rho_{m-1}$ and $\rho_m$
has no effect, since the resulting forces cancel.
Therefore, we may again assume with no loss of generality
that $\sum_{j=1}^m\rho_j=\kappa-6$.
In that case,
the image of this process under any conformal map $\psi:S\to\H$
that satisfies $\psi(w_0)\ne\infty$
is chordal-SLE$(\kappa;\rho_1,\rho_2,\dots,\rho_{m})$.
The proof is left to the dedicated reader.

\section{Associated martingales}
\label{s.mart}

In~\cite{\WernerTransformations} it was shown that a chordal SLE$(\kappa;\rho)$ process
can be viewed as an ordinary SLE$(\kappa)$ process weighted by a martingale.  In the following,
we extend this to most of the SLE-like processes discussed in
the previous sections.

Suppose that $(Y_t, t\ge 0)$ is some random process taking values in some space
$X$.  If $h:X\to [0,\infty)$ is some function such that $h(Y_1)$
is measurable and $0<\E[h(Y_1)]<\infty$, then we may weight our given probability measure
by $h(Y_1)$; that is, we may consider the probability measure $\tilde \P$ whose
Radon--Nykodim derivative with respect to $\P$ is $h(Y_1)/\E[h(Y_1)]$.
In some cases, the new law of $Y_t$ is called the Doob-transform of the
unweighted law.

In many situations one can explicitly determine the Doob-transform.
Consider, for example, a diffusion process $Y_t$ adapted to the
filtration $\mathcal F_t$ taking values
in some domain in $\R^n$.  If $h$ is as above, then
$M_t:=\Eb{h(Y_1)\md\mathcal F_t}$ is a martingale.  It turns out to be worthwhile to forget
about $h$ and consider weighting by any positive martingale.
Indeed, for every event $A\in\mathcal F_t$ we have
$\tilde \P[A]=\Eb{1_A\,M_t}/M_0$.
Girsanov's theorem (see, e.g.~\cite{\RevuzYor}) describes the behavior of
continuous martingales $N_t$ weighted by a positive continuous martingale
$M_t$ adapted to the same filtration $\mathcal F_t$.  It states that
$$
N_t-\int^t\frac{d\langle N,M\rangle_{t'}}{M_{t'}}
$$
is a local martingale under the weighted measure.
In particular, if $B_t$ is Brownian motion under $\P$ and $dM_t=a_t\,dB_t$, for some adapted
process $a_t$, then
$\tilde B_t:=B_t-\int^t (a_{t'}/M_{t'})\,dt'$ is a $\tilde\P$-Brownian motion
(note that $\langle \tilde B\rangle_t=\langle B\rangle_t$).
Thus, with respect to $\tilde P$, $B_t$ has the drift term $(a_t/M_t)\,dt$.

\begin{theorem}\label{t.mart}
Consider standard chordal $\SLE(\kappa)$.
Let $z^1,\dots,z^n$ be some collection of distinct points in $\H$.
Let $\rho_1,\dots,\rho_n\in\R$ be arbitrary.
Define $z^j_t=x^j_t+i\,y^j_t:=g_t(z^j)$.
Set
\begin{multline*}
M_t:=
\prod_{j=1}^n \Bigl(|g_t'(z^j)|^{(8-2\kappa+\rho_j)\rho_j/(8\kappa)}
(y^j_t)^{\rho_j^2/(8\kappa)}
|W_t-z_t|^{\rho_j/\kappa}
\Bigr) \, \times\\
\prod_{1\le j<j'\le n}
\bigl(\bigl|z^j_t-z^{j'}_t\bigr|
\bigl|\overline{z^j_t}-z^{j'}_t\bigr|\bigr)^{\rho_j\rho_{j'}/(4\kappa)}.
\end{multline*}
Then $M_t$ is a local martingale.
Moreover, under the measure weighted by $M_t$ (with an appropriate
stopping time) we have chordal $\SLE(\kappa;\rho_1,\dots,\rho_n)$
with force points $z^1,\dots,z^n$.
\end{theorem}

The local martingale $M_t$ was independently discovered by Marek Biskup
\cite{Biskup}.

\proof
Using \ito's formula, we may calculate
$$
dM_t =
M_t\,
\frac{1}{\kappa}\,\Re\Bigl( \sum_{j=1}^n\frac{ \rho_j}{W_t-z_t}\Bigr)
\,
dW_t
=
M_t\,
\frac{1}{\sqrt{\kappa}}\,
\Re\Bigl(
\sum_{j=1}^n\frac{ \rho_j}{W_t-z_t}\Bigr) \, dB_t
.
$$
This is laborious, but straightforward.
Thus, $M_t$ is a local martingale.
The drift term
$$
\frac{d\langle W,M\rangle_t}{M_t}
$$
from Girsanov's theorem is thus precisely the drift one has for
SLE$(\kappa;\rho_1,\dots,\rho_n)$ (since it is the product of the
coefficient of $dB_t$ in $dM_t$ and in $dW_t$ divided by $M_t$).
\QED

\begin{remark}\label{r.someboundary}
A similar variation of the theorem holds when some or all of the points $z^j$ are on the real
axis.  One only needs to replace the corresponding $y^j_t$ by $g_t'(z^j)$ in the
definition of $M_t$.
\end{remark}

We now discuss the radial version of the theorem.  Let $z^1,\dots,z^n$
be distinct points in the unit disk.  We may extend the radial maps
$g_t$ to the complement of the unit disk by Schwarz reflection,
$g_t(I_\D(z))=I_\D(g_t(z))$.  Set $z^{n+j}:=I_\D(z^j)$,
$\rho_{n+j}:=\rho_j$ and $z_t^j=g_t(z^j)$ for $j=1,\dots,2\,n$.  Then
the local martingale takes the form
\begin{equation}\label{e.Mrad}
M_t =
g'_t(0)^{q_0}
\prod_{j=1}^{2n}|W_t-z^j_t|^{\rho_j/(2\kappa)}\,|g_t'(z^j)|^{q_j}
 \prod_{1\le j<k\le 2n} |z_t^j-z_t^k|^{\rho_j\rho_k/(8\kappa)},
\end{equation}
where $q_0=(4+\bar\rho)\,\bar\rho/(8\,\kappa)$,
$\bar\rho=\sum_{j=1}^n\rho_j$ and
$q_j=(8-2\,\kappa+\rho_j)\,\rho_j/(16\,\kappa)$ for $j>0$.  Of course
$g'_t(0)=e^t$.  Likewise, if one or more of the special points $z_j$,
$j\le n$, lies on the unit circle, then the corresponding vanishing
term $|z_t^j-z_t^{n+j}|^{\rho_j^2/(8\kappa)}$ is replaced by
$|g_t'(z^j)|^{\rho_j^2/(8\kappa)}$.

\begin{remark}
In some cases these martingales are relevant to estimating the
probability of rare events in discrete models.  We mention two
examples, the first pertaining to the probability that a site is
pivotal in critical percolation, and the second pertaining to ``triple
points'' in uniform spanning trees.  In this remark we indicate,
without proof, why the martingales are relevant.

Consider critical percolation in a bounded planar domain $D$
containing the origin where $\p D$ is a simple closed curve.  Let
$z^1,\dots,z^4$ be four distinct points in clockwise order on $\p D$.
For every sufficiently small $\eps>0$, we may consider the triangular
lattice of mesh $\eps$ with a vertex at the origin and the event
$\mathcal A=\mathcal A(D,\eps,z^1,z^2,z^3,z^4)$ that the site at the
origin is pivotal for an open crossing in critical site percolation
between two opposite arcs on $\p D$ determined by the four marked
points in the union of the triangles of the grid that meet $D$.  One
could also consider the probability that the interfaces starting at an
even number $n$ of points on $\p D$ reach close to the origin.  Essentially,
the only difference is that the event $\mathcal A$ has to be defined a
bit differently, since there can be at most $3$ distinct interfaces
that meet a hexagon.

Observe that the event $\mathcal A$ is equivalent to the event that
the two (or $n/2$) percolation interfaces containing the points
$z^j$'s all visit the boundary of the hexagon dual to the site at the
origin.  As the percolation interfaces are explored, the conditional
probability of the event $\mathcal A$ is a positive martingale that
depends only on the domain cut by the explored segments of the
interfaces.  This quantity should be (nearly) conformally invariant,
modulo scaling by a power of the conformal radius of $D$ with respect
to the origin.  In particular, the martingale will not depend on
$g_t'(z^j)$ for the boundary points $z_t^j$.  Referring to
\eqref{e.Mrad}, we see that when $\kappa=6$ \cite{Smirnov}, by
selecting each $\rho_j=2$ the $g_t'(z^j)$ terms drop out, so that the
martingale $M_t$ simplifies to
$$
M_t = g'_t(0)^{(n^2-1)/12}
 \prod_{1\le j<k\le n} |z_t^j-z_t^k|^{1/3},
$$
where $W_t$ in \eqref{e.Mrad} is here written as $z^n_t$.  Letting
$f:D\to\D$ denote a conformal bijection taking the domain $D$ to the
unit disk with $f(0)=0$, it is thus reasonable to conjecture that the
probability that the $n/2$ interfaces each reach (or ``approach'' when
$n>6$) the hexagon centered at the origin will be
\begin{equation}
(\text{const}+o(1))\times (\eps f'(0))^{(n^2-1)/12} \prod_{1\le j<k\le n} |f(z^j)-f(z^k)|^{1/3}.
\end{equation}
The exponent of $(n^2-1)/12$ was shown in \cite{MR1879816}.
\old{
When the domain $D$ of interest is not the unit disk, and $v\in D$ and $f:D\to\D$ is a conformal map taking it to the unit disk, we may map $D$ to $\D$ with $g(w)=(f(w)-f(v))/(1-f(w)\overline{f(v)})$.  Then
$$g'(v)=\frac{f'(v)}{1-|f(v)|^2}\text{\ \ \ and\ \ \ }g(z^j)-g(z^k)=\frac{(f(z^j)-f(z^k))(1-|f(v)|^2)}{(1-f(z^j)\overline{f(v)})(1-f(z^k)\overline{f(v)})},$$ so the corresponding formula for the point $v$ in the domain $D$ is
$$ \text{const}\times \left|\eps f'(v)\right|^{(n^2-1)/12} \times (1-|f(v)|^2)^{(n-1)^2/12} \times \frac{\prod_{1\leq j<k\leq n}|f(z^j)-f(z^k)|^{1/3}}{\prod_{1\leq j\leq n}\left|1-f(z^j)\overline{f(v)}\right|^{(n-1)/3}}. $$
}

Consider next the uniform spanning tree with wired boundary conditions
in $D$, and let us consider the unlikely event that the tree paths
leading to $\p D$ from three different fixed neighbors of the origin
stay disjoint and hit $\p D$ at three specific points $z^1,z^2,z^3$.
After one path from a neighbor of $v$ is generated, the probability
that the next path will manage to hit its desired target is given by
the harmonic measure of that target.  In this case, rather than
looking for a martingale in which the $g_t'$ terms at the boundary
points drop out, we want a martingale in which the $g_t'$ terms at the
boundary points have an exponent of $1$.  To get this exponent of $1$
in \eqref{e.Mrad} when $\kappa=2$ \cite{\LSWlesl}, we again need to
pick the $\rho_j$'s to be $2$, giving us
$$
M_t = g'_t(0)^{(n^2-1)/4}
\prod_{1\le j< n} g_t'(z^j)
 \prod_{1\le j<k\le n} |z_t^j-z_t^k|.
$$
Again letting $f:D\to\D$ denote a conformal bijection taking the
domain $D$ to the unit disk with $f(0)=0$, this martingale suggests
the probability of this triple point (or more generally $n$-tuple point)
event is
\begin{equation}
(\text{const}+o(1))\times
(\eps f'(0))^{(n^2-1)/4}
\prod_{1\le j\le n} \eps f'(z^j)
 \prod_{1\le j<k\le n} |f(z^j)-f(z^k)|.
\end{equation}
The exponent of $(n^2-1)/4$ agrees with the value computed in \cite{MR926190,MR1757962}.
\end{remark}

\bigskip
\noindent
{\bf Acknowledgments}: We are grateful to Scott Sheffield for numerous conversations.

\bibliographystyle{nyjalpha}
\bibliography{mr,prep,notmr}

\end{document}